\documentclass[11pt]{article}
\usepackage{latexsym}
\usepackage{amssymb}
\usepackage{amsmath}
\usepackage{mathrsfs}
\usepackage{amsfonts}
\usepackage{amsthm}

\def\be{\begin{eqnarray}}
\def\ee{\end{eqnarray}}
\def\ben{\begin{eqnarray*}}
\def\een{\end{eqnarray*}}

\def\elabel#1{\label{e:#1}}

\def\sq{$\Box$}

\def\qed{\ifmmode\sq\else{\unskip\nobreak\hfil
\penalty50\hskip1em\null\nobreak\hfil\sq
\parfillskip=0pt\finalhyphendemerits=0\endgraf}\fi\par\medbreak}

\newsavebox{\junk}
\savebox{\junk}[1.6mm]{\hbox{$|\!|\!|$}}

\def\til={{\widetilde =}}

 \def\eq#1/{(\ref{#1})}

\def\eq#1/{(\ref{e:#1})}

\newcommand{\beqn}[1]{\notes{#1}%
\begin{eqnarray} \elabel{#1}}

\newcommand{\eeqn}{\end{eqnarray} }

\newcommand{\beq}[1]{\notes{#1}%
\begin{equation}\elabel{#1}}

\newcommand{\eeq}{\end{equation}}

\def\bdes{\begin{description}}
\def\edes{\end{description}}

\def\notes#1{}

\title{Some Properties of the Generalized Stuttering Poisson Distribution and its Applications}

\author{Huiming Zhang$^1${ } Lili Chu$^2${ } Yu Diao$^2$
\footnote{Supported by National University Student Innovation Test Plan. Huiming Zhang, Email:at352693585@gmail.com }\\
{\em 1.Department of Economics, Central China Normal University,}\\
{\em 152 Luoyu Street, Wuhan, 430079, China}\\
\and
{\em 2.Department of Mathematics and Statistics, Central China Normal University,}\\
{\em 152 Luoyu Street, Wuhan, 430079, China}\\
}

\date{}

\begin{document}
\bibliographystyle{plain}
\maketitle



\noindent
{\bf Abstract: }
Based on the probability generating function of stuttering Poisson distribution (SPD), this paper considers some equivalent propositions of SPD. From this, we show that some distributions in the application of non-life insurance actuarial science are SPD, such as negative binomial distribution, compound Poisson distribution etc.. By weakening condition of equivalent propositions of SPD, we define the generalized SPD. We consider cumulant estimation of  generalized SPD$'$s parameters. As an application, we use SPD with four parameters (4-th SPD) to fit auto insurance claim data. The fitting results show that 4-th SPD is more accurate than negative binomial and Poisson distribution.

\medskip

\noindent
{\bf Keywords: }stuttering Poisson distribution, probability generating function, cumulant, generalized stuttering Poisson distribution, non-life insurance actuarial science.

\medskip

\noindent
{\bf AMS 2010 Subject Classification: }
60E05, 
60E10, 
62F10, 
62P05, 


\section{Introduction}

Stuttering Poisson distribution (simply write SPD) is a non-negative discrete compound Poisson distribution (see[16],[12]), which has the feature that two or more events occur in a very short time (arrive in group or batches). For example, a man may claim for double or much compensation because he has two or more insurance policy from the insurance company. In inventory management, a customer may buy more than one goods of same kind.

{\bf Definition 1.1 } For a stochastic process $\{\xi(t),t\geq0\}$, let
\[P_n(t)=P\{\left.{\xi(t)=n}\right|\xi (0)=0\}.\]
 Similarly to some properties of Poisson process [14], stuttering Poisson process $\xi(t)$ satisfies the following properties:

{\bf Property 1:} $\xi(0)=0$.

{\bf Property 2:} $\xi(t)$ has independent increments (i.e., the numbers of events that occur in disjoint time intervals are independent) and stationary increments ( i.e., the distribution of the number of events that occur in any interval of time depends only on the length of the time interval).

{\bf Property 3:} ${P_i}(\Delta t) = \lambda {\alpha _i}\Delta t{\rm{ }} + {\rm{ }}o(\Delta t),(0 \le {\alpha _i} \le 1,i = 1,2, \cdots ,r)$, where $1\le r \le  +\infty$. The $r$ equal to finite or infinite in  the following part of this paper.

{\bf Property 3} implies that stuttering Poisson process satisfies that the probability of two or more events occur at a very short time interval is non-zero. The probability is directly proportional to the length of time interval. When $i=0$, according to the properties of independent increments , stationary increments and Chapman-Kolmogoroff equations, we have
\[{P_0}(t + \Delta t) = {P_0}(t){P_0}(\Delta t)\]
implying that
\[{P_0}(\Delta t){\rm{ = }}{{\rm{e}}^{ - \lambda \Delta t}}{\rm{ = }}1 - \lambda \Delta t + o(\Delta t),(\lambda  > 0)\]
Add terms from $P_0(\Delta t)$ to $P_r(\Delta t)$, that is
\[1 = \sum _{i = 1}^r{P_i}(\Delta t) = 1 - \lambda \Delta t + o(\Delta t) + \sum _{i = 1}^r\lambda {\alpha _i}\Delta t + {\rm{ }}o(\Delta t),\]
hence we obtain $\sum _{i = 1}^r{\alpha _i} = 1$.

{\bf Definition 1.2 } If $\xi(t)$ satisfies Property 1 ,2 and 3, then we say that $\xi(t)$ obeys the r-th stuttering Poisson distribution. Denote
\[\xi (t) \sim SP({\alpha _1}\lambda t, \cdots ,{\alpha _r}\lambda t)\]
with parameters $({\alpha _1}\lambda t, \cdots ,{\alpha _r}\lambda t) \in {\mathbb{R}^r},({\alpha _r} \ne 0)$.

Similar to Poisson law of small numbers, SPD can be deduced from the limiting distribution of multinomial distribution [19]. Let ${p_i} = \frac{{{\alpha _\iota }\lambda t}}{N}$, then the probability generating function (PGF) of stuttering Poisson distribution
$$P(s) = \mathop {\lim }\limits_{\scriptstyle N{p_i} = {\alpha _i}\lambda t \hfill \atop
  \scriptstyle   N \to \infty  \hfill} {[(1 - \sum\limits_{i = 1}^r {{p_i}} ) + \sum\limits_{i = 1}^r {{p_i}{s^i}} ]^N} = \mathop {\lim }\limits_{N \to \infty } {[1 + \frac{{\lambda t}}{N}(\sum\limits_{i = 1}^r {{\alpha _i}{s^i}}  - \sum\limits_{i = 1}^r {{\alpha _i}} )]^N} = {e^{\lambda t\sum\limits_{i = 1}^r {{\alpha _i}({s^i}}  - 1)}}.$$
where ${\alpha _i}(i = 1,2, \cdots )$ is probability density of an positive discrete distribution. With ${P_n}(t) = \frac{{{P^{(n)}}(0)}}{{n!}}$ and Fa¨¤ di Bruno fomula [15]
$$\frac{{{d^n}}}{{d{t^n}}}g[f(t)] = \sum\limits_{i = 1}^n {[\sum\limits_{\scriptstyle {k_1} +  \cdots  + {k_u} +  \cdots {k_n} = i,{\rm{ }}{k_u} \in \mathbb{N} \atop
  \scriptstyle 1 \cdot {k_1} +  \cdots  + u{k_u} +  \cdots  + n{k_n} = n} {\frac{{{\rm{ }}n!}}{{{k_1}!{k_2}! \cdots {k_n}!}}{g^{(i)}}(f(t)){{\left( {\frac{{f'(t)}}{{1!}}} \right)}^{{k_1}}}{{\left( {\frac{{f''(t)}}{{2!}}} \right)}^{{k_2}}} \cdots {{\left( {\frac{{{f^{(n)}}(t)}}{{n!}}} \right)}^{{k_n}}}]} }, $$
we have
\begin{equation}
{P_n}(t) = [{\alpha _n}\lambda t +  \cdots  + \sum\limits_{\scriptstyle {k_1} +  \cdots  + {k_u} +  \cdots {k_n} = i,{\rm{ }}{k_u} \in \mathbb{N} \atop
  \scriptstyle 1 \cdot {k_1} +  \cdots  + u{k_u} +  \cdots  + n{k_n} = n} {\frac{{\alpha _1^{{k_1}}\alpha _2^{{k_2}} \cdots \alpha _n^{{k_n}}}}{{{k_1}!{k_2}! \cdots {k_n}!}}} {(\lambda t)^i} +  \cdots  + \frac{{\alpha _1^n{{(\lambda t)}^n}}}{{n!}}]{e^{ - \lambda t}}{\rm{,}}(\lambda  > 0).
\end{equation}

Other methods to prove the expression of stuttering Poisson distribution can be obtained by the limiting distribution of multinomial distribution [19], system of differential equations [17], system of functional equation [9]. The name compound Poisson was used by W. Feller [15] and R. M. Adelson [10] to discuss distribution which PGF is ${e^{\lambda t\sum\limits_{i = 1}^r {{\alpha _i}({s^i}}  - 1)}}$. When $r=1$, SPD degenerates to Poisson distribution. When $r\geq2$, we call it non-degenerative SPD. When $r=2$, C.D. Kemp and A.W. Kemp [2] named it Hermite distribution owing to the PGF can be expanded in terms of Hermite polynomial. When $r=3(4)$, Y. C. Patel [17] said it triple (quadruple) stuttering Poisson distributions. H. P. Galliher,et al. [7] considered the demands $\xi(t)$ obey SPD with parameters
$$((1 - \alpha )\lambda ,(1 - \alpha )\alpha \lambda ,(1 - \alpha ){\alpha ^2}\lambda , \cdots )$$
of geometric distribution in inventory management theory and first name it stuttering Poisson. T.S. Moothathu, C.S. Kumar [12] considered SPD with parameters  of binomial distribution. For more application in inventory management, see [3],[6],[8]. In queue theory, A.Kuczura [1] considered requests arrive in group or batches with constant service times. R. Mitchell [11] showed  that SPD is a more exact model than the Poisson model to fit observed demand form actual historical data of several U.S. Air Force bases. D. E. Giles [4] used Hermite distribution to fit data for the number of banking and currency crises in IMF-member countries.

\section{Equivalent Propositions of SPD}
Let ${P(k, p)}$ be  a  family  of non-negative discrete distributions, which is closed  under convolution operation,  where p  (denoting  the  mean  value  of  the  distribution)  runs  over all non-negative  real numbers. L. Janossy [9] showed that ${P(k, p)}$ is SPD, The equivalent proposition of SPD is deduced from PGF and it help us to consider the generalized stuttering Poisson distribution.

{\bf Theorem 2.1 } For a discrete random variable $X=i$,$(i = 0,1,2, \cdots)$ with $P(s) = \sum _{i = 0}^r {p_i}{s^i}{\rm{ }}(\left| s \right| \le 1)$. Then, taking logarithm of PGF and expand it to a power series
\[{g_X}(s) \buildrel \Delta \over = \ln (\sum _{i = 0}^r {p_i}{s^i}) = \sum _{i = 0}^r {b_i}{s^i},\left| s \right| \le 1\]
with $\sum _{i = 0}^r {b_i} = \lambda  < \infty ,({b_i} \ge 0)$, where ${g_X}(s)$ is cumulants generation function of a discrete random variable. Then, if and only if the discrete random variable obey SPD.

{\bf Proof: } Sufficiency. With $\sum _{i = 0}^r {b_i}{s^i} = \ln (\sum _{i = 0}^r {p_i}{s^i})$ and $\sum _{i = 0}^r {b_i} = \lambda  < \infty$, we know that ${g_X}(s)$ is absolutely convergent in $\left| s \right| \le 1$. Hence $\sum _{i = 0}^r {b_i} = \ln (\sum _{i = 0}^r {p_i}) = 0$. Let ${b_i} = {a_i}\lambda $, it yield to
\[P(z) = {e^{\sum\limits_{i = 1}^r {{b_i}{s^i}}  - \sum\limits_{i = 1}^r {{b_i}} }} = {e^{\sum\limits_{i = 1}^r {{a_i}\lambda ({s^i} - 1} )}}.\]
Set $\lambda 't = \lambda $, then $X \sim SP(\lambda {\alpha _1},\lambda {\alpha _2}, \cdots ).$

Necessity. The parameters of SPD satisfy $\sum _{i = 0}^r {\alpha _i}\lambda  = \lambda  < \infty .$

{\bf Example 2.1 } Negative binomial distribution (NBD): ${P_k} = \left( {\begin{array}{*{20}{c}}{-r}\\k \\\end{array}} \right){p^r}{(p - 1)^k},(p \in (0,1),k = 1,2, \cdots)$. The PGF of NBD is ${[\frac{p}{{1 - (1 - p)s}}]^r}$. The logarithm of PGF is $ r\ln p + \sum\limits_{i = 1}^\infty  {\frac{{r{{(1 - p)}^i}}}{i}} {s^i},(\left| s \right| \le 1)$, then we have
\[r\ln p + \sum\limits_{i = 1}^r  {\frac{{r{{(1 - p)}^i}}}{i}}  < r(\ln p + \frac{{1 - p}}{p}).\]
So the NBD is equivalent to SPD with parameters $(rq,\frac{{r{q^2}}}{2},\cdots,\frac{{r{q^i}}}{i},\cdots)$. This reveals the essential properties of NBD. NBD is an important distribution in automobile insurance, which is a mixture distribution of Poisson distribution and logarithmic distribution [12].

W. Feller [15] name that the discrete random variable $X(i.i.d.,P\{ X = i\}  = {b_i},i = 0,1,2, \cdots )$ whose sums ${Y_N} = \sum _{i = 1}^N{X_i}$ with $N \sim P(\lambda t)$ is compound Poisson distribution (the accurate name should be discrete compound Poisson distribution or SPD). And PGF is
\[E(E(\left. {{s^{{Y_N}}}} \right|N = n)){\rm{ = }}E({[G(s)]^n}){\rm{ = }}{e^{\lambda [G(s) - 1]}}{\rm{ = }}{e^{\lambda \sum\limits_{i = 1}^r {{b_i}({s^i} - 1)} }}.\]

SPD have infinitely divisible property via its PGF. A nonnegative discrete distribution is called infinitely divisible
if for any $n>1$, its PGF can be represented as the n-th power of some other PGF. Thus a SPD with PGF ${e^{\lambda t\sum\limits_{i = 1}^r {{\alpha _i}({s^i}}  - 1)}}$, then this SPD can be represented as the n-th power of the other PGF ${e^{\frac{{\lambda t}}{n} \sum\limits_{i = 1}^r {{\alpha _i}({s^i}}  - 1)}}$. ${e^{\frac{{\lambda t}}{n} \sum\limits_{i = 1}^r {{\alpha _i}({s^i}}  - 1)}}$ is SPD with parameters $(\frac{{{\alpha _1}\lambda t}}{n}, \cdots ,\frac{{{\alpha _r}\lambda t}}{n})$.

Next, we obtain "compound compound Poisson" sums is also SPD in {\bf Theorem 2.2 }.

{\bf Theorem 2.2 } $X (P\{ X = i\}  = {b_i},i = 0,1,2, \cdots )$ are independent identically distributed £¨i.i.d.£© random variables. When $N\sim SP({\alpha _1}\lambda t, \cdots ,{\alpha _r}\lambda t)$, "compound compound Poisson" sums ${Y_N} = \sum _{i = 1}^N{X_i}$ is SPD.

{\bf Proof: } ${G_X}(s) = \sum _{i = 0}^r {b_i}{s^i}{\rm{,}}(\left| s \right| \le 1)$, and using double conditional expectation of ${{s^{{Y_N}}}}$
\[{P_{{Y_N}}}(s) = E({s^{{Y_N}}}) = {E_N}(E(\left. {{s^{{Y_N}}}} \right|N = n)){\rm{ = }}{E_N}({[{G_X}(s)]^n}){\rm{ = }}{e^{\lambda \sum\limits_{i = 1}^r {{\alpha _i}[G{{(s)}^i} - 1} ]}}.\]
Noticed that $\left| {{G_X}(s)} \right| \le \sum _{i = 0}^\infty {b_i} = 1$, hence
\[\lambda \sum _{i = 1}^r \left| {{\alpha _i}[{G_X}{{(s)}^i} - 1]} \right| \le \lambda \sum _{i = 1}^r {\alpha _i}{\left| {{G_X}(s)} \right|^i} + \lambda  \le 2\lambda .\]
So ${P_{{Y_N}}}(s)$ is absolutely convergent in $\left| s \right| \le 1$,
\[\lambda \sum _{i = 0}^r {\alpha _i}[G{(s)^i} - 1] = \lambda \sum _{i = 1}^r {\alpha _i}{({b_0} + {b_1}s + {b_2}{s^2} +  \cdots+  {b_r}{s^r}  )^i} - \lambda  = \sum _{i = 1}^r \lambda {c_i}{s^i} - \sum _{j = 1}^r \lambda {\alpha _j}(1 - b_0^j),\]
where $c_i$ are derived from multinomial expand. We need't to have the accurate expression of $c_i$. Let $\Sigma _{j = 1}^r {\alpha _j}(1 - b_0^j){\rm{ = }}c$, then
\[{P_{{Y_N}}}(s) = {e^{\lambda\sum\limits_{i = 1}^r {{c_i}({s^i} - 1)} }}.\]
Noticed that ${c_i} > 0$ and $\sum _{i = 1}^r \lambda {c_i} = c < \infty $, ${Y_N}$ is SPD by using {\bf Theorem 2.1 }.

In {\bf Theorem 2.2 }, a discrete compound Poisson distribution(or SPD) is a special case of discrete compound SPD (compound compound Poisson distribution) when $N \sim P(\lambda t)$. {\bf Theorem 2.2 } show that SPD is infinitely divisible distribution. It concludes that "compound $\cdots$ compound Poisson distribution" is SPD. In practice, many claims may be from superimposed events. That explains why some distributions in non-life insurance are equivalent to SPD. Besides SPD, other generalized Poisson model has wide applications in non-life insurance actuarial model and risk model such as mixed Poisson process [14],[17] and doubly stochastic Poisson processes [17].

\section{Generalized Stuttering Poisson Distribution}
In this chapter, general stuttering Poisson is defined by weakening conditions of {\bf Theorem 2.1 }.  L. Janossy [9] used independent increments, stationary increments and Chapman-Kolmogoroff equations to construction the system of functional equation
\begin{equation}
{P_i}(t + \Delta t) = \sum\limits_{k = 0}^i {{P_k}(\Delta t){P_{i - k}}(t)} {\rm{,}}(i = 0,1, \cdots )
\end{equation}
Solving (2) from one to one will deduce to (1). For example, when $i=0$, we have ${P_0}(t + \Delta t) = {P_0}(t){P_0}(\Delta t)$, the solution of ${P_0}(t)$ is ${P_0}(t){\rm{ = }}{{\rm{e}}^{ - \lambda t}}$.

When $i=1$, we have ${P_1}(t + \Delta t) = {P_0}(t){P_1}(\Delta t) + {P_1}(t){P_0}(\Delta t)$, the solution is ${P_1}(t) = {\alpha _1}t{e^{ - \lambda t}}$.

When $i=2,\cdots$, by the system of functional equation (2) we have
\begin{equation}
{P_n}(t) = [{\alpha _n}\lambda t +  \cdots  + \sum\limits_{\scriptstyle {k_1} +  \cdots  + {k_u} +  \cdots {k_n} = i,{\rm{ }}{k_u} \in \mathbb{N} \atop
  \scriptstyle 1 \cdot {k_1} +  \cdots  + u{k_u} +  \cdots  + n{k_n} = n} {\frac{{\alpha _1^{{k_1}}\alpha _2^{{k_2}} \cdots \alpha _n^{{k_n}}}}{{{k_1}!{k_2}! \cdots {k_n}!}}} {(\lambda t)^i} +  \cdots  + \frac{{\alpha _1^n{{(\lambda t)}^n}}}{{n!}}]{e^{ - \lambda t}}{\rm{,}}(\lambda  > 0),
\end{equation}
where $  (\lambda  \ge 0,{\alpha _i}\in\mathbb{R},n =0, 1,2, \cdots ,r)$. Janossy [9] proved (3) by mathematical induction.

There are no nonnegative restriction in (3), implies
\[{P_i}(\Delta t) = \lambda {\alpha _i}\Delta t{\rm{ }} + {\rm{ }}o(\Delta t),(i = 1,2, \cdots)\]
we obtain $\sum _{i = 1}^\infty{\alpha _i} = 1$.

Since ${\alpha _i}$ are not necessarily to be nonnegative, now we suppose that ${\alpha _i}$ may take negative value and satisfy $\sum _{i = 1}^\infty {\alpha _i} = 1$ to have a new distribution family.

{\bf Definition 3.1 } Generalized stuttering Poisson distribution (GSPD): For a discrete random variable $X (X=i,i=0,1,2, \cdots )$ , the form of PGF is ${e^{\lambda \Sigma _{i = 1}^\infty ({\alpha _i}{s^i} - 1)}}$ and satisfies $\lambda  > 0,\sum\limits_{i = 1}^\infty  {{\alpha _i}}  = 1,\sum\limits_{i = 1}^\infty  {\left| {{\alpha _i}} \right|}  < \infty $. When ${\alpha _i} \equiv 0$, if ${\alpha _i} \ge r + 1$. We name it r-th generalized stuttering Poisson distribution (GSPD).

It is obvious that SPD is a subfamily of GSPD.

{\bf Example 3.1 } {\bf Theorem 3.1 }  show that the probability of zero occurrence is more than positive occurrences. For example, the cumulant of Bernoulli distribution $P(s) = p + (1 - p)s{\rm{, }}(p > 0.5)$ is
\[\ln [p + (1 - p)s] =  - \ln p + \frac{{1 - p}}{p}s - \frac{1}{2}{(\frac{{1 - p}}{p})^2}{s^2}{\rm{ + }}\frac{1}{3}{(\frac{{1 - p}}{p})^3}{s^3}{\rm{ + }} \cdots .\]

\section{Statistic of Generalized Stuttering Poisson Distribution}
Cumulants $\kappa_n$,moments $m_n$ and central moments $c_n$ of GSPD are deduced from probability generating function and moment generating function. It also can use in SPD. A discrete random variable $X (P\{ X = i\}  = {b_i},i = 0,1,2, \cdots )$ with PGF $P(s) = \sum _{i = 0}^\infty {p_i}{s^i},(\left| s \right| < 1)$, where moment generating function is ${M_X}(s) = P({e^s})$. Expanding ${e^{sX}}$ with Taylor series at zero, we have
\[{M_X}(s) = E(\sum\limits_{n = 0}^\infty  {\frac{{{{(sX)}^n}}}{{n!}}} ) = \sum\limits_n^\infty  {\frac{{E{X^n}}}{{n!}}} {s^n} \buildrel \Delta \over = \sum\limits_{n = 0}^\infty  {\frac{{{m_k}}}{{n!}}} {s^n},(\left| s \right| \le 1)\]

{\bf Definition 4.1 }  Cumulants generating function of a random variable is
\[{g_X}(s) \buildrel \Delta \over = \ln ({M_X}(s)) = \sum\limits_{n = 0}^\infty  {{\kappa _n}\frac{{{s^n}}}{{n!}}} ,(\left| s \right| \le 1)\]
where coefficients ${\kappa _n}\,(n = 0,1,2, \cdots )$ is n-th cumulants. It is explicit that ${\kappa _0} = 0$.

{\bf Theorem 4.1 } $\xi (t) \sim GSP({\alpha _1}\lambda t, \cdots ,{\alpha _r}\lambda t)$, then the n-th cumulant of $\xi (t)$ is
\begin{align}
{\kappa _n} = \sum _{i = 0}^r{\alpha _i}\lambda t{i^n}.
\end{align}
{\bf Proof } Cumulants generating function of  GSPD is ${g_X}(s)=\sum _{i = 0}^r{\alpha _i}\lambda t({e^{is}} - 1)$, Expanding ${e^{sX}}$ with Taylor series at zero, that is
\begin{align*}
\sum\limits_{i = 1}^r {{\alpha _i}\lambda t({e^{is}} - 1)} & =  - \lambda t + \lambda t[{\alpha _1}(1 + \sum\limits_{j = 1}^\infty  {\frac{{{s^j}}}{{j!}}} ) + {\alpha _2}(1 + \sum\limits_{j = 1}^\infty  {\frac{{{{(2s)}^j}}}{{j!}}} ) +  \cdots  + {\alpha _k}(1 + \sum\limits_{j = 1}^\infty  {\frac{{{{(ks)}^j}}}{{j!}}} )]\\
& = (\sum\limits_{i = 1}^r {{\alpha _i}\lambda ti} )s + (\sum\limits_{i = 1}^r {{\alpha _i}\lambda t{i^2}} )\frac{{{s^2}}}{{2!}} +  \cdots  + (\sum\limits_{i = 1}^r {{\alpha _i}\lambda t{i^l}} )\frac{{{s^l}}}{{l!}} +  \cdots.
\end{align*}
Comparing coefficients of $t^n$, we have ${\kappa _n} = \sum _{i = 0}^r{\alpha _i}\lambda t{i^n}$. When $r \to \infty $, $\kappa_n$ may divergence.

{\bf Theorem 4.2 } If $\xi (t) \sim GSP({\alpha _1}\lambda t, \cdots ,{\alpha _r}\lambda t)$, then the recursion formula of nth moments $m_n$ is
\begin{align}
{m_{n + 1}} = \sum\limits_{j = 0}^n {\left( {\begin{array}{*{20}{c}}
   n  \\
   r  \\
\end{array}} \right){\kappa _{n + 1 - j}}{m_j}}  = \sum\limits_{j = 0}^n {\left( {\begin{array}{*{20}{c}}
   n  \\
   r  \\
\end{array}} \right)(\sum\limits_{i = 1}^r {{\alpha _i}\lambda t {i^{n + 1 - j}}} ){m_j}} ,{\rm{ (}}{\kappa _1} = {m_1} = \sum\limits_{i = 1}^r {{\alpha _i}\lambda  t i} {\rm{ ) }}.
\end{align}

{\bf Proof } Expanding $\ln [{M_X}(s)]$ with Taylor series at zero, that is
\begin{align*}
\ln [{M_X}(s)] & = \ln (1 + \sum\limits_{n = 1}^\infty  {\frac{{{m_n}}}{{n!}}} {s^n}) = \sum\limits_{n = 1}^\infty  {\frac{{{m_n}}}{{n!}}} {s^n} - \frac{1}{2}{(\sum\limits_{n = 1}^\infty  {\frac{{{m_n}}}{{n!}}} {s^n})^2} +  \cdots  + \frac{{{{( - 1)}^{i - 1}}}}{i}{(\sum\limits_{n = 1}^\infty  {\frac{{{m_n}}}{{n!}}} {s^n})^i} +  \cdots\\
& = {m_1}s + \frac{{{m_2} - m_1^2}}{{2!}}{s^2} + \frac{{{m_3} - 3{m_1}{m_2} + 2m_1^3}}{{3!}}{s^3} + \frac{{{m_4} - 4{m_3}{m_1} - 3m_2^2 + 12{m_2}m_1^2 - 6m_1^4}}{{4!}}{s^4} +  \cdots.
\end{align*}
By the definition of cumulant, we have
$${\kappa _1} = {m_1} = EX = \sum _{i = 1}^r{\alpha _i}\lambda ti{\rm{}} \Rightarrow {m_1} = \sum _{i = 1}^r{\alpha _i}\lambda ti,$$
$${\kappa _2} = {m_2} - m_1^2 = E{(X - EX)^2} = \sum _{i = 1}^r{\alpha _i}\lambda t{i^2}{\rm{}} \Rightarrow {m_2} = \sum _{i = 1}^r{\alpha _i}\lambda t{i^2} + {(\sum _{i = 1}^r{\alpha _i}\lambda ti)^2},$$
$${\kappa _3} = {m_3} - 3{m_1}{m_2} + 2m_1^3 = E{(X - EX)^3}{\rm{ = }}\sum _{i = 1}^r{\alpha _i}\lambda t{i^3}{\rm{  }} \Rightarrow {m_3} = {\kappa _3} + 3{m_1}{m_2} - 2m_1^3,$$
$${\kappa _4} = {m_4} - 4{m_1}{m_3} - 3m_2^2 + 12m_1^2{m_2} - 6m_1^4 = E{(X - EX)^4} - 3{[E{(X - EX)^2}]^2}, \cdots$$
By taking the derivative of both side of ${M_X}(s){\rm{ = }}{{\rm{e}}^{{g_X}(s)}}$ again and again with respect to $s$, using Leibniz formula we obtain
\begin{align}
M_X^{(n + 1)}(s) = \sum\limits_{i = 0}^n {\left( {\begin{array}{*{20}{c}}
   n  \\
   i  \\
\end{array}} \right){{[g_X^{(1)}(s)]}^{(n - i)}}M_X^{(i)}(s)}.
\end{align}
Substitute $s=0$ to the n-th derivatives of ${g_X}(s)$ and ${M_X}(s)$, hence ${\left. {M_X^{(n)}(s)} \right|_{s = 0}} = {m_n}$  and ${\left. {g_X^{(n)}(t)} \right|_{s = 0}} = {\kappa _n}$. Substituting cumulants and moments into (6), we have
\begin{align}
{m_{n + 1}} = \sum\limits_{j = 1}^n {\left( {\begin{array}{*{20}{c}}
   n  \\
   j  \\
\end{array}} \right){\kappa _{n + 1 - j}}{m_j}}  = \sum _{j = 1}^n\left( {\begin{array}{*{20}{c}}
   n  \\
   i  \\
\end{array}} \right)(\sum _{i = 1}^r{\alpha _i}\lambda t {i^{n + 1 - j}}){m_j}.
\end{align}
The relationship between cumulants and moment is
\begin{align}
{\kappa _n} & = {m_n} - \sum\limits_{j = 1}^{n - 1} {\left( {\begin{array}{*{20}{c}}
   {n - 1}  \\
   i  \\
\end{array}} \right){\kappa _{n - j}}{m_j}} {\rm{, (}}{\kappa _1} = {m_1}{\rm{)}}\\
& \buildrel \Delta \over =f(m_1,m_2,\cdots,m_n).
\end{align}

{\bf Remark} Higher cumulants $(n\geq4)$ are different to central moment when $n$ is more than 4. Arguing from ${\kappa _n} = \frac{{{d^n}}}{{d{t^n}}}\ln [{M_X}(s)]$ and Fa¨¤ di Bruno formula, it can deduce to ${\kappa _n}$ too.

{\bf Example 4.2 } An alternative approach to compute the mean and variance of Bernoulli distribution
\[{\kappa _n} = \sum\limits_{i = 1}^\infty  {{\alpha _i}\lambda {i^n}}  = \sum\limits_{i = 1}^\infty  {\frac{{{{(1 - p)}^i}}}{{ - i{p^i}\ln p}}( - \ln p){i^n}}  = \sum\limits_{i = 1}^\infty  {{{(\frac{{1 - p}}{p})}^i}} {i^{n - 1}}\]
\[ \Rightarrow EX{\rm{ = }}\sum\limits_{i = 1}^\infty  {{{(\frac{{1 - p}}{p})}^i}} {\rm{ = }}1 - p,DX = \sum\limits_{i = 1}^\infty  {i{{(\frac{{1 - p}}{p})}^i}}  = p - {p^2}.\]

{\bf Theorem 4.3 } If $\xi (t) \sim GSP({\alpha _1}\lambda t, \cdots ,{\alpha _r}\lambda t)$, then the recursion formula of nth moments $c_n$ is
\begin{align}
{c_{n + 1}} = \sum\limits_{j = 0}^n {\left( {\begin{array}{*{20}{c}}
   n  \\
   i  \\
\end{array}} \right)\kappa _{n + 1 - j}^*{c_j}} ,{\rm{ }}({c_0} = 1,{c_1} = 1{\rm{ }},\kappa _n^* = \left\{ {\begin{array}{*{20}{c}}
   {{\kappa _n},{\rm{ }}n \ne 1}  \\
   {0,{\rm{   }}n = 1}  \\
\end{array}} \right.)
\end{align}

{\bf Proof } Expanding $\ln [{M_{X - EX}}(s)]$ with Taylor series at zero and comparing the coefficients of cumulant generating function, we obtain
\[\kappa _1^* = 0,{\rm{ }}\kappa _2^* = \sum _{i = 1}^r{\alpha _i}\lambda t{i^2},{\rm{ }}\kappa _3^* = \sum _{i = 1}^r{\alpha _i}\lambda t{i^3},{\rm{ }}\kappa _4^* = \sum _{i = 1}^r{\alpha _i}\lambda t{i^4} + 3{(\sum _{i = 1}^r{\alpha _i}\lambda t{i^2})^2}, \cdots \]
Since
$$M_{X - EX}^{(n + 1)}(s) = \sum\limits_{i = 0}^n {\left( {\begin{array}{*{20}{c}}
   n  \\
   i  \\
\end{array}} \right){{[g_{X - EX}^{(1)}(s)]}^{(n - i)}}M_{X - EX}^{(i)}(s)}, $$
substitute $s=0$ to the n-th derivatives of ${g_{X - EX}}(s)$ and ${M_{X - EX}}(s)$, hence ${\left. {M_{X - EX}^{(n)}(s)} \right|_{s = 0}} = {c_n},{\left. {g_{X - EX}^{(n)}(s)} \right|_{s = 0}} = \kappa _n^*$. Similarly to the proof in {\bf Theorem 4.2 }, displacing ${\kappa _n}$ with $\kappa _n^*$, we have (10).

\section{Cumulant Estimation of Generalized Stuttering Poisson Distribution}
Y. C. Patel [18] gave moment estimator of the parameters of Hermite distribution. Y.C.Patel [19] estimates the parameters of the trip and quadruple stuttering Poisson distributions with maximum likelihood estimation moment estimation, and mixed moment estimation of the parameter. Use Sameple moments $\hat m_1$ and central moments $\hat c_2$ and $\hat c_3$, we have
\begin{align}
\left( {\begin{array}{*{20}{c}}
   {{{\hat m}_1}}  \\
   {{{\hat c}_2}}  \\
   {{{\hat c}_3}}  \\
\end{array}} \right) \buildrel \Delta \over = \left( {\begin{array}{*{20}{c}}
   1 & 2 & 3  \\
   1 & 4 & 9  \\
   1 & 8 & {27}  \\
\end{array}} \right)\left( {\begin{array}{*{20}{c}}
   {{{\hat \alpha }_1}}  \\
   {{{\hat \alpha }_2}}  \\
   {{{\hat \alpha }_3}}  \\
\end{array}} \right)\lambda t{\rm{ }} \Rightarrow {\rm{ }}\left( {\begin{array}{*{20}{c}}
   {{{\hat \alpha }_1}}  \\
   {{{\hat \alpha }_2}}  \\
   {{{\hat \alpha }_3}}  \\
\end{array}} \right)\lambda t = \left( {\begin{array}{*{20}{c}}
   3 & {{\textstyle{{ - 5} \over 2}}} & {{\textstyle{1 \over 2}}}  \\
   {{\textstyle{{ - 3} \over 2}}} & 2 & {{\textstyle{{ - 1} \over 2}}}  \\
   {{\textstyle{1 \over 3}}} & {{\textstyle{{ - 1} \over 2}}} & {{\textstyle{1 \over 6}}}  \\
\end{array}} \right)\left( {\begin{array}{*{20}{c}}
   {{{\hat m}_1}}  \\
   {{{\hat c}_2}}  \\
   {{{\hat c}_3}}  \\
\end{array}} \right),
\end{align}
where $E{\hat \alpha _i} = {\alpha _i} + O({\textstyle{1 \over n}})\,,(i = 1,2,3)$.

When $n\geq4$, From the computing in Theorem 4.2, Higher cumulants are different to central moment or moment. Central moment or moment are nonlinear combination of ${\alpha _i}(i = 1,2, \cdots )$. Thus it is difficult to estimate the parameters by using central moment or moment estimation. {\bf Theorem 4.2 } implies that $\kappa_n$  is linear combination of ${\alpha _i}(i = 1,2, \cdots )$.

Therefore, firstly, we use sample moment $\kappa_n$ to calculate in (4). From (4), when $r < \infty $,
$\kappa_n$ is convergence. Secondly, by solving the following system of linear equations of ${\hat \alpha _i}(i = 1,2, \cdots ,n)$ by means of 0-th to n-th cumulants formula
\begin{align}
\left( {\begin{array}{*{20}{c}}
   {{{\hat \kappa }_0}}  \\
   {{{\hat \kappa }_1}}  \\
    \vdots   \\
   {{{\hat \kappa }_n}}  \\
\end{array}} \right) \buildrel \Delta \over = \left( {\begin{array}{*{20}{c}}
   1 & 1 &  \cdots  & 1  \\
   0 & 1 &  \cdots  & n  \\
    \vdots  &  \vdots  &  \ddots  &  \vdots   \\
   0 & {{1^n}} &  \cdots  & {{n^n}}  \\
\end{array}} \right)\left( {\begin{array}{*{20}{c}}
   { - 1}  \\
   {{{\hat \alpha }_1}}  \\
    \vdots   \\
   {{{\hat \alpha }_n}}  \\
\end{array}} \right){\rm{ }},
\end{align}
where the coefficients matrix in (12) is invertible Vandermonde matrix. By solving the linear system equation we have
\begin{align}
{\hat \alpha _i} = \sum\limits_{j = 1}^n {{b_{ij}}} {\hat \kappa _j}, (i = 1,2, \cdots ,n)
\end{align}

Assuming the samples to the power of $n$, $\xi _i^n(i = 1,2, \cdots ,m)$ are i.i.d.. Let samples n-th moment ${A_n} = \frac{1}{m}\sum\limits_{i = 1}^m {\xi _i^n} $, arguing from Khintchine's law of large numbers, we obtain $\mathop {\lim }\limits_{n \to \infty } P\{ \left| {{A_n} - {m_n}} \right| < \varepsilon \}  = 1$  for all $\epsilon>0$.
From relationship between cumulants and moment (4), for all $\epsilon>0$. ${f(x_1,x_2,\cdots,x_n)}$ is continuous of several variables, it implies
\begin{align}
\mathop {\lim }\limits_{n \to \infty } P\{ \left| {{\kappa_n} - {{\hat \kappa _n}}} \right| < \varepsilon \}  =\mathop {\lim }\limits_{n \to \infty } P\{ \left| {{f(m_1,m_2,\cdots,m_n)} - {f(A_1,A_2,\cdots,A_n)}} \right| < \varepsilon \}=1
\end{align}
Thus we prove cumulant estimation is consistent estimate by the linear relation in (13).

\section{ Applications}
R. M. Adelson [17] put forward the recursion formula of SPD's probability density function by using Leibniz formula
\begin{align}
{P_{j + 1}}(t) = \frac{1}{{j + 1}}[{\alpha _1}\lambda {P_j}(t) + 2{\alpha _2}\lambda {P_{j - 1}}(t) +  \cdots  + (j + 1){\alpha _{j + 1}}\lambda {P_0}(t)],{\rm{ }}{P_0}(t) = {e^{ - \lambda t}},
\end{align}
(15) avoid tediously computing the sum of much index in (3) by recursion relation. There are no nonnegative restriction to ${\alpha _i}(i = 1,2, \cdots )$, so (15) can be used in GSPD.

Now we use SPD fitting according to auto  insurance claims data (the car insurance claims data of the following table 1 from [17]), and then we compare the goodness of fit with some other distributions

From data in Table 1, Total insurance policies are $n=106974$. The probability of zero claim policies is far greater than 0.5. Obviously it is zero-inflated data. The number of the insurance policy of i-th is ${x_j}(j = 1,2, \cdots ,106974)$, so the mean value and the 2-th and the 3-th central moment of the insurance policy claims rate is
\[{m_1} = 0.1010806364,{\rm{ }}{c_2} = 0.1074468102,{\rm{ }}{c_3} = 0.1216468798.\]
According to the (8) and (11),  we have
\[{\hat \alpha _1} = {\rm{0}}{\rm{.97255}},{\hat \alpha _2} = {\rm{0}}{\rm{.02496}},{\hat \alpha _3} = {\rm{0}}{\rm{.00249}}.\]
Thus we can infer that the probability of claims of customer who buy two copies of the same insurance is only $2.496{\rm{\% }}$, and three copies of the same insurance is only $0.249{\rm{\% }}$. Employing recursion relation (15) and (8), we obtain ${\hat p_i}$. And then figure out $n{\hat p_i}$. Analogously, consider quadruple SPD fitting, in this case we have
\[{\tilde \alpha _1} = {\rm{0}}{\rm{.97151}}, {\tilde \alpha _2} = {\rm{0}}{\rm{.02703}}, {\tilde \alpha _3} = {\rm{0}}{\rm{.00112, }}{\tilde \alpha _4}{\rm{ = 0}}{\rm{.00034}}.\]

In Table 1, we assume the data come from Poisson distribution, triple SPD (by recursion formula (15)), quadruple SPD (by recursion formula (15)), negative binomial distribution, respectively, and then estimate the probability of different numbers of claims.
\begin{table}
\caption{The comparison of auto insurance claims data of different distributions fitting effect (moment estimation or cumulant estimation)} 
\centering          
\begin{tabular}{c c c c c c c c}    
\hline\hline                        
$i$ & Methods & 0 & 1 & 2 & 3 & 4& $i>4$  \\ [0.5ex]  
\hline                      
  $v_i$ & {\bf Observed frequency } & {\bf 96978 } & {\bf 9240 } & {\bf 704 } & {\bf 43 } & {\bf 9 } & {\bf 0 } \\
  $n{\hat p_i}$ & Estimate by Poisson & 96689.5 & 9773.5 & 494.5 & 16.6  & 0.4 & 0 \\
  $n{\hat p_i}$ & Estimate by triple SPD & 96974.1 & 9256.0 & 679.2 & 60.4 & 4.0 & 0.3 \\
  $n{\hat p_i}$ & Estimate by quadruple SPD & 96977.3 & 9243.2 & 697.6 & 49.1 & 6.1 & 0.7 \\
  $n{\hat p_i}$ & Estimate by NPD & 96985.4 & 9222.5 & 711.7 & 50.7 & 3.5 & 0.2 \\[1ex]        
\hline          
\end{tabular}
\label{table:nonlin}    
\end{table}

Constructing test statistical:$\eta  = \sum\limits_{i = 0}^4 {\frac{{v_i^2}}{{n{{\hat p}_i}}} - n} $, by calculating, we get
\[{\eta _{PD}} = {\rm{345}}{\rm{.1250}},{\rm{ }}{\eta _{3-SPD}} = {\rm{12}}{\rm{.5786}},{\rm{ }}{\eta _{4-SPD}} = {\rm{2}}{\rm{.8963}},{\rm{ }}{\eta _{NBD}} = {\rm{10}}{\rm{.1294}}.\]

From Pearson's chi-squared test theory, in one hand, $\chi _4^2(0.01) = 12.277$, given significant level of 0.1 we accept that claims data obey quadruple SPD or negative binomial distribution. In the other hand, $\chi _4^2(0.5) = 3.357$, given significant level of 0.5 we accept that claims data obey quadruple SPD fitting. It is thus clear that quadruple SPD fitting effect is better than that of negative binomial distribution fitting effect, and NBD is better than triple SPD. The goodness of Poisson distribution model is worst in those four distributions.

\section*{Acknowledgments}
The authors want to thank Chengming Sun, Prof. Hui Zhao and Prof. Yinbao Chen for their support and helpful comments.



\end{document}